\documentclass[11pt]{article}
\usepackage[all]{xy}
\usepackage{amsfonts}
\usepackage{amsthm}
\usepackage{enumerate}
\usepackage{graphicx}
\usepackage{mathrsfs}
\usepackage{bm}
\usepackage{cite}
\usepackage[colorlinks,linkcolor=blue,citecolor=blue]{hyperref}
\usepackage{amssymb,amsmath} % font used for R in Real numbers
\usepackage{hyperref}
\usepackage{makecell}
\usepackage{tikz}
\usepackage{pgf}
\usepackage{tikz}
\usetikzlibrary{patterns}
\usepackage{pgffor}
\usepackage{pgfcalendar}
\usepackage{pgfpages}
\usepackage{shuffle,yfonts}
\usepackage{mathtools}
\DeclareFontFamily{U}{shuffle}{}
\DeclareFontShape{U}{shuffle}{m}{n}{ <-8>shuffle7 <8->shuffle10}{}

 \allowdisplaybreaks
\usetikzlibrary{arrows,shapes,chains}
 \allowdisplaybreaks

\catcode`!=11
\let\!int\int \def\int{\displaystyle\!int}
\let\!lim\lim \def\lim{\displaystyle\!lim}
\let\!sum\sum \def\sum{\displaystyle\!sum}
\let\!sup\sup \def\sup{\displaystyle\!sup}
\let\!inf\inf \def\inf{\displaystyle\!inf}
\let\!cap\cap \def\cap{\displaystyle\!cap}
\let\!max\max \def\max{\displaystyle\!max}
\let\!min\min \def\min{\displaystyle\!min}
\let\!frac\frac \def\frac{\displaystyle\!frac}
\catcode`!=12

\let\oldsection\section
\renewcommand\section{\setcounter{equation}{0}\oldsection}

\theoremstyle{plain}

\newtheorem{theorem}{Theorem}[section]
\newtheorem{conj}[theorem]{Conjecture}
\newtheorem{proposition}[theorem]{Proposition}
\newtheorem{lemma}[theorem]{Lemma}
\newtheorem{cor}[theorem]{Corollary}
\theoremstyle{remark}

\newcommand{\bbar}[1]{\overline{#1}}

\newcommand{\dep}{\operatorname{dep}}
\newcommand{\wt}{\operatorname{wt}}
\begin{document}
%%%%%%%%%%%%%%%%%%%% title %%%%%%%%%%%%%%%%%%%%%%%%%%%%%%%%%%%%%%%%%%%%%%%%
\title{\bf An Analytic Evaluation of Euler Sharp Sums}
\author{
{Zhenlu Wang$^{a,}$\thanks{Email: zhenluwang@tzc.edu.cn},\quad Ce Xu$^{b,}$\thanks{Email: cexu2020@ahnu.edu.cn}\quad{and}\quad Jianqiang Zhao${}^{c,}$\thanks{Email: zhaoj@ihes.fr}}\\[1mm]
a. \small School of Artificial Intelligence, Taizhou University, \\ \small Taizhou 318000, {Zhejiang,} P.R. China\\
b. \small School of Mathematics and Statistics, Anhui Normal University,\\ \small Wuhu 241002, {Anhui,} P.R. China\\
c. \small Department of Mathematics, The Bishop's School, La Jolla, CA 92037, USA}

\date{}
\maketitle

\noindent{\bf Abstract.} We give a direct analytic evaluation of the Euler sharp sums $\zeta^\sharp(2,\{1\}^{d-3},\overline{2})$ for all even integers $d\ge4$, expressing them in terms of ordinary multiple zeta values of depth at most two. The proof reduces these sharp sums to a harmonic-number series and evaluates it using a generating-function identity, thereby proving the even-depth case of a conjecture of Xu and Zhao.
\medskip

\medskip


\section{Introduction}

A finite sequence of positive integers is called an \emph{index}. The weight and depth of the index $\boldsymbol{k}=(k_1,\ldots,k_r)$ are
$\wt(\boldsymbol{k}):=k_1+\cdots+k_r$ and $\dep(\boldsymbol{k}):=r$,
respectively.

For an index $\boldsymbol{k}=(k_1,\ldots,k_r)$ with $k_r\ge2$, the \emph{multiple zeta value} (MZV) and the \emph{multiple zeta-star value} (MZSV) are defined by
\begin{align*}%\label{eq:mzv}
 \zeta(k_1,\ldots,k_r):=\sum_{1\le n_1<\cdots<n_r}
 \frac1{n_1^{k_1}\cdots n_r^{k_r}}\quad\text{and}\quad \zeta^\star(k_1,\ldots,k_r):=\sum_{1\le n_1\le\cdots\le n_r}
 \frac1{n_1^{k_1}\cdots n_r^{k_r}},
\end{align*}
respectively. The systematic study of MZVs began with the works of Hoffman \cite{H1992} and Zagier \cite{DZ1994} in the early 1990s. Due to their surprising and sometimes mysterious appearance in the study of many branches of mathematics and theoretical physics, these special values have attracted a lot of attention and interest in the past three decades.  A natural level-two generalization of MZVs is provided by \emph{Euler sums} (also called \emph{alternating multiple zeta values}), which arise by inserting factors $(-1)^{n_j}$. These sums attracted considerable attention in the 1990s through their appearances in the evaluation of Feynman integrals. Since then, their rich algebraic relations and arithmetic properties have motivated extensive research; see \cite{Zhao2016} for details.

Following Xu and Zhao \cite{XuZhao2026}, for an index $\boldsymbol{k}=(k_1,\ldots,k_r)$ and $\boldsymbol{\varepsilon}=(\varepsilon_1,\ldots,\varepsilon_r)\in\{\pm1\}^r$ with $(k_r,\varepsilon_r)\ne(1,1)$, the \emph{Euler sum} is defined by
$$\zeta(\boldsymbol{k};\boldsymbol{\varepsilon})=\zeta(k_1,\ldots,k_r;\varepsilon_1,\ldots,\varepsilon_r):=\sum_{0<n_1<\cdots<n_r}\frac{\varepsilon_1^{n_1}\cdots\varepsilon_r^{n_r}}{n_1^{k_1}\cdots n_r^{k_r}}.$$
Obviously, ordinary MZVs are recovered when all $\varepsilon_j=1$. As usual, we write $\overline{k_j}$ in place of $k_j$ when $\varepsilon_j=-1$, and leave $k_j$ unbarred when $\varepsilon_j=1$. For example, $\zeta(1,\overline{1})=\zeta(1,1;1,-1)$.

Let $x_0,x_1,x_{-1}$ be the one forms defined by
$$x_0:=\frac{dt}{t},\qquad x_1:=\frac{dt}{1-t},\qquad
 x_{-1}:=\frac{dt}{-1-t}.$$
Then Euler sums admit the iterated integral expressions
\begin{align}\label{eq:ES-integral}
 \zeta(k_1,\ldots,k_r;\varepsilon_1,\ldots,\varepsilon_r)=\int_0^1 x_{\eta_1}x_0^{k_1-1}\cdots x_{\eta_r}x_0^{k_r-1},
\end{align}
where $\eta_j=\varepsilon_j\cdots\varepsilon_r$ for $1\le j\le r$. The theory of iterated integrals was developed firstly by K.T. Chen in the 1960's \cite{KTChen1971,KTChen1977}. It has played important roles in the study of algebraic topology and algebraic geometry in the past half century. Its simplest form is
$$\int_{a}^b f_p(t)dtf_{p-1}(t)dt\cdots f_1(t)dt:=\int_{a<t_p<\cdots<t_1<b}f_p(t_p)f_{p-1}(t_{p-1})\cdots f_1(t_1)dt_1dt_2\cdots dt_p.$$
As with ordinary MZVs, Euler sums satisfy two fundamental product relations: the shuffle product arising from their iterated integral representations and the stuffle (or harmonic) product arising from their series definitions. For example, let $\overline{\mathbb N}:=\{\overline n:n\in\mathbb N\}$. For any $a,b\in\mathbb N\cup\overline{\mathbb N}$, the stuffle product gives
$$\zeta(a)\zeta(b)=\zeta(a,b)+\zeta(b,a)+\zeta(a\oplus b),$$
where for positive integers $m,n$, $\oplus$ is defined by
\begin{align*}
m\oplus n=\overline m\oplus\overline n:=m+n,\quad
m\oplus\overline n=\overline m\oplus n:=\overline{m+n}.
\end{align*}

For all $k_1,\ldots,k_r\in\mathbb N\cup\overline{\mathbb N}$, the \emph{Euler sharp sum} is defined by
\begin{align}\label{eq:ESS-def}
\zeta^\sharp(k_1,\ldots,k_r):=\sum_{\boldsymbol{q}=(k_1\circ\cdots\circ k_r)}2^{\operatorname{dep}(\boldsymbol{q})}\zeta(\boldsymbol{q}),
\end{align}
where each $\circ$ is chosen to be either a comma ``,'' or ``$\oplus$'', and $\dep(\boldsymbol q)$ denotes the depth of $\boldsymbol q$. Thus, Euler sharp sums can be regarded as weighted linear combinations of Euler sums obtained by merging adjacent components in all possible ways. For example,
$$\zeta^\sharp(2,\overline 3)=4\zeta(2,\overline 3)+2\zeta(\overline 5).$$
Euler sharp sums are closely related to interpolated Euler sums at $t=1/2$, whose MZV counterparts were introduced by Yamamoto \cite{Yamamoto2013}. More precisely,
$$\zeta^\sharp(\boldsymbol{k})=2^{\operatorname{dep}(\boldsymbol{k})}\zeta^{1/2}(\boldsymbol{k}).$$

Euler sharp sums arise naturally in the study of identities relating multiple zeta-star values (MZSVs) to Euler sums. A related identity obtained by Linebarger and Zhao \cite{LinebargerZhao2015} is
\begin{align}\label{eq:LZ-identity}
\zeta^\sharp(\overline{2},\{1\}^{d-3},\overline{2})=\zeta^\star(1,d),\qquad d\ge3,
\end{align}
where $\{a\}^n$ denotes $n$ repetitions of $a$. Moreover, Zhao \cite[Theorem~1.4]{Zhao2016TwoOne} established another family of two-one type identities, one special case of which gives
$$\zeta^\sharp(\{1\}^{n-2},\overline{2})=-\zeta(n),
\qquad n\ge2.$$
These identities illustrate how certain Euler sharp sums of arbitrarily large depth can be reduced to ordinary MZVs or MZSVs.

Recently, Xu and Zhao \cite{XuZhao2026} conjectured an explicit evaluation of the family of Euler sharp sums $\zeta^\sharp(2,\{1\}^{d-3},\overline{2})$, which differs from the Euler sharp sum on the left-hand side of \eqref{eq:LZ-identity} only in that its first entry is unbarred rather than barred. Such an evaluation gives an explicit MZV expression (called \emph{unramified}) for a height-one family of multiple mixed values (MMVs)\cite{XuZhao2022}; the corresponding motivic MMVs were shown in \cite{XuZhao2026} to be unramified for even $d$. In the present paper, we prove the conjectured identity for all even integers $d\ge4$. Our main result is stated as follows.

\begin{theorem}\label{thm:main}
For any even integer $d\ge4$, we have
\begin{align}\label{eq:main}
 \zeta^\sharp(2,\{1\}^{d-3},\bbar2)
 =\frac12\zeta(d+1)
 -2\sum_{\substack{m+n=d+1\\2\mid n}}\zeta(m,n).
\end{align}
\end{theorem}
The proof first reduces the Euler sharp sum to a single harmonic-number series, valid for all integers $d\ge3$, and then evaluates this series for even $d\ge4$ by means of a generating-function identity.

Furthermore, based on our numerical computations using GP-Pair, we find that the following family of Euler sharp sums is also unramified. Thus, we propose the following conjecture:
\begin{conj}
Let $b$ be odd and $a,c \le 2$. Then the Euler sharp sum  \(\zeta^\sharp (\{2\}^a, \{1\}^b, \{\bar{2}\}^c)\) is unramified.
\end{conj}
Clearly, when \(a=c=1\), the conjecture reduces to Theorem \ref{thm:main}. However, the methods presented in this paper appear to fail for the remaining cases of the conjecture.

\section{Proof of the main theorem}\label{sec:proof}
We begin with an iterated-integral formula for the sharp sum. By the integral expression of Euler sums \eqref{eq:ES-integral} and the definition of Euler sharp sums \eqref{eq:ESS-def}, we have
\begin{align}\label{eq:sharpint}
 \zeta^\sharp(k_1,\ldots,k_r)=2\int_0^1
 x_{\eta_1}x_0^{|k_1|-1}
 \prod_{j=2}^r(x_0+2x_{\eta_j})x_0^{|k_j|-1},
\end{align}
where $|k_j|$ denotes the underlying positive integer if $k_j$ is barred.
Indeed, replacing a separator by $x_0$ merges the adjacent components,
whereas choosing $2x_{\eta_j}$ retains it. A coarsening of depth $h$
therefore appears with coefficient $2\cdot2^{h-1}=2^h$.

Let $$H_n:=\sum_{j=1}^n\frac{1}{j}$$ be the \emph{classical harmonic numbers} and $H_0:=0$. We also use the \emph{Pochhammer symbol}, defined by $(a)_0:=1$ and $(a)_n:=a(a+1)\cdots(a+n-1)$ for $n\ge1$. The following series expansion will be used in the proof of the main theorem.
\begin{lemma}\label{lem:key-expansion}
For $|x|<1$, we have
\begin{align}
\frac{1}{\sqrt{1-x}}\log\frac{2}{1+\sqrt{1-x}}=\sum_{n=1}^{\infty} c_n\bigl(H_{2n}-H_n\bigr)x^n,
\label{eq:key-expansion}
\end{align}
where $$c_n=\frac{1}{4^n}\binom{2n}{n}.$$
\end{lemma}

\begin{proof}
Differentiating the binomial series
$(1-x)^{-\alpha}=\sum_{n\ge0}(\alpha)_nx^n/n!$
with respect to $\alpha$, then taking $\alpha=1/2$, gives
\begin{align*}
-\frac{\log(1-x)}{\sqrt{1-x}}
=\sum_{n\ge1}c_n(2H_{2n}-H_n)x^n.
\end{align*}
On the other hand, the integral representation
$H_n=\int_0^1(1-t^n)/(1-t)\,dt$ gives, for $0<x<1$,
\begin{align*}
\sum_{n\ge1}c_nH_nx^n
&=\int_0^1\frac{1}{1-t}
 \left(\frac{1}{\sqrt{1-x}}
 -\frac{1}{\sqrt{1-tx}}\right)dt\\
&=\frac{2}{\sqrt{1-x}}
 \log\frac{1+\sqrt{1-x}}{2\sqrt{1-x}}.
\end{align*}
Here the second equality follows by putting
$y=\sqrt{1-tx}$; nonnegative summands justify the interchange
of sum and integral. Both sides extend analytically to $|x|<1$,
with the branches normalized at $x=0$.
The assertion now follows by taking half the difference of the
two displayed generating functions and using
$\log(1-x)=2\log\sqrt{1-x}$.
\end{proof}

Then, we reduce the Euler sharp sum to a harmonic number series.
\begin{proposition}\label{prop:harmonic}
For any integer $d\ge3$, we have
\begin{align}\label{eq:harmonic}
 \zeta^\sharp(2,\{1\}^{d-3},\bbar2)
 =-2\sum_{n\ge1}\frac{H_{2n}-H_n}{n^d}.
\end{align}
\end{proposition}
\begin{proof}
By \eqref{eq:sharpint}, we have
\begin{align*}%\label{eq:qint}
 \zeta^\sharp(2,\{1\}^{d-3},\bbar2)=2\int_0^1x_{-1}x_0\gamma^{d-2}x_0,
\end{align*}
where $$\gamma=x_0+2x_{-1}=\frac{dt}{t}-2\frac{dt}{1+t}.$$
Fix the variables $u<v$ belonging to the two occurrences of $x_0$.
For a single one-form $\omega$, symmetry of the integration simplex gives
$$\int_{u<t_1<\cdots<t_s<v}\prod_{j=1}^s\omega(t_j)
 =\frac1{s!}\left(\int_u^v\omega\right)^s.$$
Applying this identity to $\gamma$, and noting that
$$\int_0^u x_{-1}=-\log(1+u),\qquad\int_u^v\gamma=\log\frac{v(1+u)^2}{u(1+v)^2},$$
we obtain
\begin{align*}%\label{eq:uv}
\zeta^\sharp(2,\{1\}^{d-3},\bbar2)=-\frac2{(d-2)!}\int_{0<u<v<1}\frac{\log(1+u)}{uv}\left(\log\frac{v(1+u)^2}{u(1+v)^2}\right)^{d-2} du\,dv.
\end{align*}
The quadratic change of variables
\begin{align*}
x=\frac{4u}{(1+u)^2},\qquad
y=\frac{4v}{(1+v)^2}
\end{align*}
sends $0<u<v<1$ onto $0<x<y<1$ and satisfies
$$\frac{du}{u}=\frac{dx}{x\sqrt{1-x}},\quad
 \frac{dv}{v}=\frac{dy}{y\sqrt{1-y}},\quad
 \frac{v(1+u)^2}{u(1+v)^2}=\frac yx.$$
Consequently, we obtain
\begin{align}\label{eq:xy}
 \zeta^\sharp(2,\{1\}^{d-3},\bbar2)=-\frac2{(d-2)!}\int_{0<x<y<1}
 \frac{\log\!\bigl(2/(1+\sqrt{1-x})\bigr)}
 {xy\sqrt{1-x}\sqrt{1-y}}
 \left(\log\frac yx\right)^{d-2}dx\,dy.
\end{align}
Substituting \eqref{eq:key-expansion} and
$(1-y)^{-1/2}=\sum_{b\ge0}c_by^b$ in \eqref{eq:xy} gives
\begin{align*}
\zeta^\sharp(2,\{1\}^{d-3},\bbar2)=-\frac2{(d-2)!}\sum_{a\ge1}\sum_{b\ge0}
 c_ac_b(H_{2a}-H_a)
 \int_{0<x<y<1}x^{a-1}y^{b-1}
 \left(\log\frac yx\right)^{d-2}dx\,dy.
\end{align*}
For $a\ge1$ and $b\ge0$, the substitution $x=yt$ gives
\begin{align*}
\int_{0<x<y<1} x^{a-1}y^{b-1}\left(\log\frac{y}{x}\right)^{d-2}\,dx\,dy&=\int_0^1y^{a+b-1}dy\int_0^1t^{a-1}(-\log t)^{d-2}dt\notag\\
&=\frac{1}{a+b}\frac{(d-2)!}{a^{d-1}}.
%\label{eq:basic-double-integral}
\end{align*}
Moreover, by the binomial expansion and the Euler beta integral, we have
\begin{align*}%\label{eq:beta}
 \sum_{b\ge0}\frac{c_b}{a+b}
 =\int_0^1\frac{t^{a-1}}{\sqrt{1-t}}\,dt
 =B\left(a,\frac12\right),
\end{align*}
where the \emph{Beta function} \(B(x,y) = \int_0^1 t^{x-1}(1-t)^{y-1} dt\) for \(\Re(x), \Re(y) > 0\) is related to the Gamma function by \(B(x,y) = \frac{\Gamma(x)\Gamma(y)}{\Gamma(x+y)}\).

Since
\begin{align*}
 c_a B\left(a,\frac12\right)=\frac1a,
\end{align*}
it follows that
\begin{align*}
\zeta^\sharp(2,\{1\}^{d-3},\overline2)=-2\sum_{a\ge1}\frac{c_a(H_{2a}-H_a)}{a^{d-1}}\sum_{b\ge0}\frac{c_b}{a+b}=-2\sum_{a\ge1}\frac{H_{2a}-H_a}{a^d},
\end{align*}
which proves \eqref{eq:harmonic}.
\end{proof}

For $|z|<1$, we define the following generating function
\begin{align*}
 S(z)&=\sum_{1\le k<\ell}\frac{k}{(k^2-z^2)(\ell^2-z^2)}.
\end{align*}
Let $\ell=k+r$. Then we have
\begin{align*}
S(z)
&=\sum_{k,r\ge1}
\frac{k}{(k^2-z^2)((k+r)^2-z^2)}\\
&=\sum_{k,r\ge1}
\frac{k}{r(2k+r)}
\left(
\frac1{k^2-z^2}
-
\frac1{(k+r)^2-z^2}
\right).
\end{align*}
Since
\begin{align*}
\sum_{r=1}^{\infty}\frac{k}{r(2k+r)}
=\frac12\sum_{r=1}^{\infty}
\left(\frac1r-\frac1{r+2k}\right)
=\frac12H_{2k},
\end{align*}
and
\begin{align*}
\sum_{k=1}^{n-1}\frac{k}{(n-k)(n+k)}
&=\frac12\sum_{k=1}^{n-1}
\left(\frac1{n-k}-\frac1{n+k}\right)=H_n-\frac12H_{2n}-\frac1{4n},
\end{align*}
we obtain
\begin{align}\label{eq:partial}
 S(z)=C(z)-A(z)+\frac14Z(z),
\end{align}
where
\begin{align*}
 A(z)=\sum_{n\ge1}\frac{H_n}{n^2-z^2},\quad
 C(z)=\sum_{n\ge1}\frac{H_{2n}}{n^2-z^2},\quad
 Z(z)=\sum_{n\ge1}\frac1{n(n^2-z^2)}.
\end{align*}

 We can now complete the proof of Theorem \ref{thm:main}.
\begin{proof}[{\bf Proof of Theorem~\ref{thm:main}}]
By \eqref{eq:partial}, for $|z|<1$, we have
\begin{align*}
\sum_{q\ge1}\left(-2\sum_{n\ge1}\frac{H_{2n}-H_n}{n^{2q}}\right)z^{2q-2}=-2(C(z)-A(z))=\frac12Z(z)-2S(z).
\end{align*}
On the other hand, we have
\begin{align*}
Z(z)&=\sum_{q\ge1}\zeta(2q+1)z^{2q-2},
\end{align*}
and
\begin{align*}
S(z)&=\sum_{q\ge1}\left(\sum_{\substack{m+n=2q+1\\2\mid n}}\zeta(m,n)\right)z^{2q-2},
\end{align*}
Comparing the coefficients of $z^{2q-2}$ for $q\ge2$, and applying \eqref{eq:harmonic} with $d=2q$, we obtain
\begin{align*}
\zeta^\sharp(2,\{1\}^{2q-3},\overline2)=\frac12\zeta(2q+1)-2\sum_{\substack{m+n=2q+1\\2\mid n}}
\zeta(m,n),
\end{align*}
as desired.
\end{proof}

According to \cite[Theorem 2.2]{Flajolet-Salvy}, \(\zeta(1,d)\) (for \(d>1\)) can be expressed as a combination of products of Riemann zeta values, and \cite[Theorem 7.2]{Flajolet-Salvy} shows that \(\zeta(1,\overline{2q})\) can also be expressed as a combination of products of Riemann zeta values. Moreover, a direct computation from Proposition~\ref{prop:harmonic} yields
\[
\zeta^\sharp(2,\{1\}^{d-3},\bbar2) = (2-2^d)\zeta(1,d) - 2^d\zeta(1,\bar{d}) + \zeta(d+1).
\]
Therefore, we have the following corollary.
\begin{cor}
For any even integer $d\ge4$, the \(\zeta^\sharp(2,\{1\}^{d-3},\bbar2)\) can be expressed as a combination of products of Riemann zeta values.
\end{cor}

{\bf Declaration of competing interest.}
The authors declare that they have no known competing financial interests related to this paper.

{\bf Data availability.}
No data were used for the research described in the article.

{\bf Acknowledgments.} C. Xu is supported by the General Program of Natural Science Foundation of Anhui Province (Grant No. 2508085MA014). J. Zhao was supported by the Jacobs Prize from the Bishop's School.

\end{document}